\documentclass{amsart}
\usepackage{amscd}
\usepackage{amsmath}
\usepackage{enumerate}
\usepackage{cite}
\usepackage[dvips]{graphicx}
\usepackage{amssymb,amsxtra,amstext,amsbsy,amsthm,latexsym,calligra,mathrsfs}    
\usepackage{amsfonts,euscript,delarray,enumitem}
\usepackage{fancyhdr,nameref,amscd}
\usepackage[all,cmtip]{xy}
\usepackage[letter,center]{crop}
\usepackage{geometry}
\usepackage{tikz,tikz-cd}
\usetikzlibrary{arrows.meta} 
\tikzcdset{every label/.append style = {font = \small}}
\usetikzlibrary{matrix,arrows}

\newtheorem{theorem}{Theorem}[section]
\newtheorem{lemma}[theorem]{Lemma}
\newtheorem{proposition}[theorem]{Proposition}
\newtheorem{corollary}[theorem]{Corollary}
\newtheorem{definition}[theorem]{Definition}

\newtheorem{remark}[theorem]{Remark}

\newtheorem{example}[theorem]{Example}

 \newtheorem{theorem*}{Theorem}
 \newtheorem{corollary*}[theorem*]{Corollary}
 \newtheorem{proposition*}[theorem*]{Proposition}

\def\proof{\removelastskip\par\medskip \noindent {\sc Proof.}\enspace}
\def\endproof{\hbox{ }{\qed}\hfill\par\medskip}

\newcommand{\CC}{{\mathbb{C}}}

\newcommand{\set}[2]{\ensuremath{\left \{\,{#1}\mid {#2}\,\right\}}}
\newcommand{\smvee}{\raise0.9ex\hbox{$\scriptscriptstyle\vee$}}

\newcommand{\Ss}[1]{\mathcal{O}_{#1}}
\newcommand{\Cs}[1]{\omega_{#1}}
\newcommand{\Rs}{\tilde{X}}

\newcommand{\Sf}[1]{\mathcal{#1}}
\newcommand{\dimc}[1]{\dim_{\mathbb{C}}\left(#1\right)}

\DeclareMathOperator{\im}{Im}

\DeclareMathOperator\coker{coker}
\DeclareMathOperator{\Exts}{\mathscr{E}\text{\kern -3pt {\calligra\large xt}}\,}
\DeclareMathOperator{\Homs}{\mathscr{H}\text{\kern -3pt {\calligra\large om}}\,}

\DeclareMathOperator{\Ext}{\text{Ext}\,}

\begin{document}
\title[On the blow-up of a normal singularity at MCM modules]{On the blow-up of a normal singularity at maximal Cohen-Macaulay modules}

\author{Agust\'in Romano-Vel\'azquez}
\address{School of Mathematics, Tata Institute of Fundamental Research, Homi Bhabha Road, Mumbai 400 005, India}
\email{agustin@math.tifr.res.in}

\thanks{The author is partially supported by ERCEA 615655 NMST Consolidator Grant, CONACYT CB 2016-1 Num. 286447, CONACYT 253506, FORDECYT 265667 and by TIFR Visiting Fellow.}

\subjclass[2010]{Primary: 13C14, 13H10, 14E16, 32S25, 32S05}
\begin{abstract}
Raynaud and Gruson developed the theory of blowing-up an algebraic variety $X$ along a coherent sheaf $M$ in the sense that there exists a blow-up $X'$ of $X$ such that the ``strict transform" of $M$ is flat over $X'$ and the blow-up satisfies an universal (minimality) property. However, not much is known about the singularities of the blow-up. In this article, we prove that if $X$ is a normal surface singularity and $M$ is a reflexive $\mathcal{O}_{X}$-module, then such a blow-up arises naturally from the theory of McKay correspondence. We show that the normalization of the blow-up of Raynaud and Gruson is obtained by a resolution of $X$ such that the full sheaf $\mathcal{M}$ associated to $M$ (i.e., the reflexive hull of the pull-back of $M$) is globally generated and then contracting all the components of the exceptional divisor not intersecting the first Chern class of $\mathcal{M}$. Moreover, we prove that if $X$ is Gorenstein and $M$ is special in the sense of Wunram and Riemenschneider (generalized in a previous work by Bobadilla and the author), then the blow-up of Raynaud and Gruson is normal. Finally, we use the theory of matrix factorization  developed by Eisenbud, to give concrete examples of such blow-ups.
\end{abstract}
\maketitle

\section{Introduction} 
Let $X$ be an algebraic variety and $M$ be a coherent $\Ss{X}$-module. The ground field is always assumed to be the complex numbers. The Raynaud-Gruson flattening theorem~\cite{Ray,Ray2} states that under some hypothesis, there exists a finitely presented closed subscheme of $X$ (depending on $M$) such that the blow-up $f\colon X' \to X$ along the subscheme satisfies the property: $f^* M / \mathrm{tor}$ is flat over $X'$. The construction is universal in the sense that for any morphism $\eta \colon Y \to X$ for which the pull-back of $M$ to $Y$ is flat modulo torsion, the morphism $\eta$ factors through $X'$. The existence of such blow-ups satisfying the universal minimality condition has numerous applications. Raynaud~\cite{Ray} uses this blow-up to prove the Chow's lemma (i.e., given $X$ separated over $S$, we can find a blow-up $X'$ of $X$ such that $X'$ is quasi-projective over $S$). Hironaka~\cite{Hironaka} recovers this blow-up and the Chow's lemma in the context of analytic geometry. Abramovich, Karu, Matsuki and Włodarczyk~\cite{Abra} use the Chow's lemma to prove the weak factorization conjecture for birational maps (i.e., a birational map between complete non-singular varieties is a composition of blow-ups and blow-downs along smooth centers). Campana~\cite{Campana} uses the flattening theorem to study the geometry, arithmetic, and classification of compact Kähler manifolds. Hassett and Hyeon~\cite{Hassett} use the flattening theorem to study the log canonical models for the moduli space of curves. Rossi~\cite{Ros}  and later Villamayor~\cite{Villa} investigate the case of the blow-up at a coherent module $M$ over a (possible non-reduced) ring $R$ with applications to the flattening of projective morphisms and related problems in singularity theory like the Nash transformation (the closure of the graph of the Gauss map which assigns each regular point to its tangent space considered as an element in a Grassmanian manifold).

In general, the Raynaud-Gruson blow-up may be very difficult to describe. Moreover, the blow-up is not regular in general, with little known about the singularities. In this article we prove that this blow-up has a very nice description in the case of $M$ a maximal Cohen-Macaulay module and $X$ a normal surface singularity. For this, we generalize some techniques and ideas given in~\cite{BoRo} together with a careful study of the Raynaud-Gruson blow-up. In this setting, given $M$ a maximal Cohen-Macaulay $\Ss{X}$-module, by~\cite{BoRo} there exists an unique resolution $\pi \colon \Rs \to X$ called the \emph{minimal adapted resolution} associated to $M$ such that the associated \emph{full sheaf} $\Sf{M}:=(\pi^*M)^{\smvee \smvee}$ is generated by global sections and $\Rs$ is the minimal resolution of $X$ satisfying this property, where $\left(-\right)^{\smvee}$ denotes the dual with respect to the structure sheaf. Recall, full sheaf was first defined by Esnault~\cite{Es} for rational singularities and generalized by Kahn~\cite{Ka} for normal surface singularities. We prove (Theorem~\ref{Theo:BlowUpMMinAdapt}):

\begin{theorem*}
\label{Theo:IntroBlowUpMMinAdapt}
Let $(X,x)$ be the germ of a normal surface singularity. Let $M$ be a reflexive  $\Ss{X}$-module of rank $r$. Let $\pi\colon \Rs \to X$ be the minimal adapted resolution associated to $M$ with exceptional divisor $E$. Let $E_1,\dots,E_n$ be the irreducible components of $E$ and $\Sf{M}:=(\pi^*M)^{\smvee \smvee}$ be the full sheaf associated to $M$. Then, the normalization of the Raynaud-Gruson blow-up of $X$ at $M$ is obtained by contracting the irreducible components $E_i$ such that $c_1\left(\Sf{M} \right) \cdot E_j =0$.
\end{theorem*}

In the case of normal, Gorenstein surface singularities and maximal Cohen-Macaulay modules we prove that the Raynaud-Gruson blow-up is related to the McKay correspondence. The McKay correspondence was constructed by McKay~\cite{McK} and a conceptual geometric understanding of the correspondence was achieved by a series of papers by Gonzalez-Springberg and Verdier~\cite{GoVe}, Artin and Verdier~\cite{AV} and by Esnault and Kn\"orrer~\cite{EsKn}. This correspondence gives a bijection between the isomorphism classes of non-trivial indecomposable reflexive modules and the irreducible components of the exceptional divisor of the minimal resolution of a rational double point. Later Wunram~\cite{Wu} generalized the McKay correspondence to any rational singularity using the notion of \emph{specialty}. If the singularity is Gorenstein we generalize in~\cite{BoRo} the definition of specialty as follows: a maximal Cohen-Macaulay module $M$ of rank $r$ over a normal surface singularity $X$ is called \emph{special} if the minimal adapted resolution $\pi \colon \Rs \to X$ has the property that the dimension over $\CC$ of the module $R^1 \pi_* \left(\pi^*M\right)^{\smvee}$ is equal to $rp_g$, where $p_g$ is the \emph{geometric genus} of $X$. Using the notion of specialty we prove the following result (Theorem~\ref{Th:BlMXNormal}):

\begin{theorem*}
\label{Th:IntroBlMXNormal}
Let $(X,x)$ be the germ of a normal Gorenstein surface singularity and $M$  be a  special module. Let $f \colon \mathrm{Bl}_M(X) \to X$ the Raynaud-Gruson blow-up of $X$ at the module $M$. Then, $\mathrm{Bl}_M(X)$ is normal.
\end{theorem*}

Theorem~\ref{Theo:IntroBlowUpMMinAdapt} and Theorem~\ref{Th:IntroBlMXNormal} generalize the results of Gustavsen and Ile~\cite{Gus} to the case of normal surface singularities. They prove that the blow-up at any maximal Cohen-Macaulay module in a rational surface singularity is a partial resolution dominated by the minimal resolution. This follows by the following two properties of rational singularities:
\begin{itemize}
    \item the minimal resolution of a rational singularity is the minimal adapted resolution of every maximal Cohen-Macaulay module,
    \item the blow-up of a \emph{complete ideal} (in the sense of Lipman~\cite{Li}) is always normal.
\end{itemize}
 For a general normal surface singularity both assertions fail. Indeed, by~\cite{BoRo} we know that in general in a Gorenstein, normal, surface singularity, its minimal resolution is not the minimal adapted resolution of every maximal Cohen-Macaulay module. Moreover, in Example~\ref{ex:Ex2} we provide a normal, Gorenstein, surface singularity and a maximal Cohen-Macaulay module such that the Raynaud-Gruson blow-up is not normal. In order to prove both theorems we use the properties of the minimal adapted resolution and the techniques developed in~\cite{BoRo}. 

As a consequence of Theorem~\ref{Theo:IntroBlowUpMMinAdapt} and  Theorem~\ref{Th:IntroBlMXNormal} we generalize the McKay correspondence given by Artin and Verdier~\cite{AV} as follows (Corollary~\ref{Cor:PrincipalNuevo}):

\begin{corollary*}
\label{Cor:IntroPrincipalNuevo}
Let $(X,x)$ be a normal Gorenstein surface singularity. Then, there exists a bijection between the following sets: 
\begin{enumerate}
\item The set of special, indecomposable $\Ss{X}$-modules up to isomorphism.
\item The set of irreducible divisors $E$ over $x$, such at any resolution of $X$ where $E$ appears, the Gorenstein form has neither zeros nor poles along $E$.
\item The set of partial resolutions $\psi\colon Y \to X$ with irreducible exceptional divisor $E$ such that:
    \begin{enumerate}
        \item the Gorenstein form has neither zeros nor poles along $E \setminus \text{Sing $Y$}$.
        \item the partial resolution is dominated by a resolution such that the Gorenstein form has neither zeros nor poles along its exceptional divisor.
    \end{enumerate}
\end{enumerate}
\end{corollary*}

In the last section we use the theory of matrix factorizations developed by Eisenbud~\cite{Ei} to compute some explicit examples of blow-up at reflexive modules. In particular we prove that the \emph{fundamental module} in the hypersurface singularity given by $f=x^3+z^3+y^3$ is not a special module.

The organization of this paper is as follows: In \S~\ref{Sec:Pre} we give preliminary results about full sheaves over normal, Gorenstein surface singularities, the Raynaud-Gruson blow-up in coherent sheaves and complete ideals. In \S~\ref{Sec:Adap} we define the notion of adapted resolutions over any normal singularity of arbitrary dimension. In \S~\ref{sec:Blowup} we prove our main results and generalize some results given in~\cite{BoRo}. In \S~\ref{Sec:Matrix} we recall some basics on matrix factorization and compute explicit examples of the Raynaud-Gruson blow-up of maximal Cohen-Macaulay modules.

\section{Preliminaries}
\label{Sec:Pre}
In this section we recall basics on full sheaves over Gorenstein singularities, the Raynaud-Gruson blow-up in coherent sheaves and complete ideals. We assume basic familiarity with dualizing sheaves, modules and normal surface singularities, see~\cite{BrHe,Har1,Ne,Ishi} for more details.

\subsection{Setting and notation} Throughout this article, we denote by $(X,x)$ either a complex analytic normal surface germ, or the spectrum of
a normal complete $\mathbb{C}$-algebra of dimension $2$. In few instances it will denote the spectrum of a normal complete $\mathbb{C}$-algebra of dimension $n$.

In this situation $X$ has a \emph{dualizing sheaf} $\omega_X$, and we also denote by $\omega_X$ its stalk at $x\in X$, which is called the \emph{dualizing module} of the ring $\Ss{X,x}$ (see~\cite[Chapter~5~\S~3]{Ishi} for more details). If $(X,x)$ is a Gorenstein normal singularity, then the dualizing module coincides with $\Ss{X,x}$. 
Let
\begin{equation*}
    \pi \colon \Rs\to X,
\end{equation*}
be a resolution of singularities. The exceptional divisor is denoted by $E:=\pi^{-1}(x)$, with irreducible components $E_1,\dots,E_m$.

If $(X,x)$ is a Gorenstein surface singularity, there is a $2$-form $\Omega_{\Rs}$ which is meromorphic in $\Rs$, and has neither zeros nor poles in $\Rs\setminus E$; this form is called the \emph{Gorenstein form}. Let $\text{div}(\Omega_{\Rs})=\sum q_iE_i$ be the divisor associated with the Gorenstein form. 
The coefficients $q_i$ are independent of the choice of the form $\Omega_{\Rs}$ with these properties.
\begin{definition}
\label{def:smallresgor}
Let  $\pi \colon \Rs\to X$ be a resolution of a normal Gorenstein surface singularity. The {\em canonical cycle} is defined as $Z_k:=\sum_i -q_iE_i$, where the $q_i$ are the coefficients defined above. 

We say that $\Rs$ is {\em small with respect to the Gorenstein form} if $Z_k$ is greater than or equal to $0$. 

The {\em geometric genus} of $X$ is defined to be the dimension as a $\mathbb{C}$-vector space of $R^1\pi_*\Ss{\Rs}$ for any resolution.
\end{definition}

Following Reid~\cite{Reid} we have the following definition:

\begin{definition}
Let $(X,x)$ be the germ of a normal surface singularity. A \emph{partial resolution} is a proper birational morphism $\pi\colon Y\to X$ with $Y$ a normal variety.
\end{definition}

\subsection{Cohen-Macaulay modules and reflexive modules}
Let $X$ be a normal variety. Let $\Homs_{\Ss{X}}(\bullet,\bullet)$ and $\Exts^i_{\Ss{X}}(\bullet,\bullet)$ be the sheaf theoretic Hom and Ext functors. 
The dual of an $\Ss{X}$-module $M$ is denoted by $M^{\smvee}:=\Homs_{\Ss{X}}(M,\Ss{X})$ and its $\omega_X$-dual is $\Homs_{\Ss{X}}(M,\omega_X)$. An $\Ss{X}$-module $M$ is called \emph{reflexive} (resp. $\omega_X$-\emph{reflexive}) if the natural homomorphism from $M$ to $M^{\smvee \smvee}$ (resp. to $\Homs_{\Ss{X}}(\Homs_{\Ss{X}}(M,\omega_X),\omega_X)$)  is an isomorphism. 

An $\Ss{X}$-module $M$ is called \emph{Cohen-Macaulay} if for every $y \in X$ the depth of the stalk $M_y$ is equal to the dimension of the support of the module $M_y$. If the depth of $M_y$ is equal to the dimension of $\Ss{X,y}$, 
then the module $M_y$ is called \emph{maximal Cohen-Macaulay}. A module is {\em indecomposable} if it cannot be written as a direct sum of two non-trivial submodules. 

\subsection{Full sheaves and the minimal adapted resolution} Let $(X,x)$ be the germ of a normal surface singularity and
\begin{equation*}
 \pi \colon \Rs \to X,   
\end{equation*}
be a resolution. Recall, the following definition of full sheaves as in~\cite[Definition~1.1]{Ka}.
\begin{definition}
An $\Ss{\Rs}$-module $\Sf{M}$ is called \emph{full} if there is a reflexive $\Ss{X}$-module $M$ such that 
$\Sf{M} \cong \left(\pi^* M\right)^{\smvee \smvee}$. We call $\Sf{M}$ the full sheaf associated to $M$.
\end{definition}

Another important notion is the concept of specialty. Wunram~\cite{Wu} and Riemenschneider~\cite{Rie} defined a special full sheaf as a full sheaf for which its dual has the first cohomology group equal to zero. In \cite{BoRo} the author and Bobadilla generalized this definition as follows:
\begin{definition}\label{def:especial}
\label{def:espmodule}
Let $M$ be a reflexive $\Ss{X}$-module of rank $r$ and $\Sf{M}$ be the full sheaf associated to $M$. The full sheaf $\Sf{M}$ is called \emph{special} if $\dimc{R^1 \pi_* \left(\Sf{M}^{\smvee}\right)} = rp_g$.
 
We say that $M$ is \emph{a special module} if for any resolution, the full sheaf associated to $M$ is special.
\end{definition}

Let $M$ be a reflexive $\Ss{X}$-module. The minimal adapted resolution associated to $M$ was defined in \cite{BoRo} and it plays a crucial role for the classification of special reflexive modules. Recall,

\begin{definition}
\label{def:minadap}
Let $M$ be a reflexive $\Ss{X}$-module. The minimal resolution $\pi \colon \Rs\to X$ for which the associated full sheaf $\left(\pi^* M\right)^{\smvee \smvee}$ is generated by global 
sections is called the {\em minimal adapted resolution} associated to $M$.  
\end{definition}

\subsection{Blowing up at coherent sheaves}
In this paper we use the description of the  Raynaud-Gruson blow-up given by Villamor~\cite{Villa}. Let $R$ be a domain with quotient field $K$. The fractional ideals are the finitely generated $R$-submodules of $K$. Two fractional ideals $J_1$ and $J_2$ are isomorphic if and only if there exists some $k$ in $K$ different from zero such that $J_1=kJ_2$. The norm of a module is defined in the class of all fractional ideals modulo isomorphism as follows:

\begin{definition}[{\cite[p.~123]{Villa}}]
Let $M$ be a finitely generated $R$-module of rank $r$. The \emph{norm} of $M$ is the class
\begin{equation*}
    \| M \|_R := \im\left( \bigwedge^r M \to \bigwedge^r M \otimes K \cong K
    \right)/\sim,
\end{equation*}
where $\sim$ denotes the isomorphism as fractional ideals.
\end{definition}

The blow-up of $R$ at the module $M$ is described as follows:

\begin{theorem}[{\cite[Theorem~3.3]{Villa}}]
\label{Th:BlowUpM}
Let $M$ be a finitely generated $R$-module of rank $r$. There exists a blow-up:
\begin{equation*}
    f\colon \mathrm{Bl}_M(R) \to \text{Spec}(R),
\end{equation*}
with the following properties:
\begin{enumerate}
    \item The sheaf $f^*M / \mathrm{tor}$ is a locally free sheaf of $\Ss{\mathrm{Bl}_M(R)}$-modules of rank $r$.
    \item (Universal property) For any morphism $\sigma \colon Z \to \text{Spec}(R)$ such that $\sigma^* M / \mathrm{tor}$ is a locally free sheaf of $\Ss{Z}$-modules of rank $r$, there exists an unique morphism $\beta\colon Z \to \mathrm{Bl}_M(R)$ such that $f \circ{} \beta = \sigma$.
\end{enumerate}
\end{theorem}
As mentioned before, Villamayor constructed the above blow-up under more general assumptions, but in our setting Theorem~\ref{Th:BlowUpM} is sufficient. Furthermore, in our situation $R$ is a domain, hence this blow-up is a proper birational morphism (see \cite{Villa} and \cite{Ros} for details). 
\begin{remark}
\label{remark:BlowupM}
Recall that, the blow-up stated in Theorem~\ref{Th:BlowUpM} was constructed by Villamayor by taking the blow-up at the fractional ideal $\|M\|_R$ (at any representative). This fact will be used later in the article.
\end{remark}

\subsection{Complete ideals} In this section we recall some basic notions of complete ideals. See~\cite{Li} for more details. Let $X$ be an integral scheme with sheaf of rational functions $\mathcal{R}_X$. Let $\mathcal{J}$ be a quasi-coherent $\Ss{X}$-submodule of $\mathcal{R}_X$. Denote by 
\begin{equation*}
\mathcal{A}:= \bigoplus_{n \geq 0} \mathcal{J}^n,\, (\mathcal{J}^0 := \Ss{X}), \quad \text{and} \quad  \mathcal{R}:=\bigoplus_{n \geq 0} \mathcal{R}_X^n, \, (\mathcal{R}_X^0 := \mathcal{R}_X),    
\end{equation*}
where $\mathcal{J}^n$ (resp. $\mathcal{R}_X^n$) is the product (as fractional ideal) of $n$ copies of $\mathcal{J}$ (resp. $\mathcal{R}_X$). Hence, the sheaves $\mathcal{A}$ and $\mathcal{R}$ are quasi-coherent graded $\Ss{X}$-algebras. Let $\mathcal{A}'$ be the integral closure of $\mathcal{A}$ in $\mathcal{R}$. Denote by $\mathcal{J}_n$ the image of the following composition:
\begin{equation*}
    \mathcal{A}' \subset \mathcal{R} \xrightarrow{\mathrm{pr}_n} \mathcal{R}_X,
\end{equation*}
where $\mathrm{pr}_n$ is the $n$-th projection.

\begin{definition}
We say that $\mathcal{J}_1$ is the \emph{completion} of $\mathcal{J}$. Furthermore, $\mathcal{J}$ is complete if $\mathcal{J}=\mathcal{J}_1$.
\end{definition}

The following two lemmas will play a crucial role later.

\begin{lemma}[{\cite[Lemma~5.2]{Li}}]
\label{lema:blowupcomplete}
Let $X$ be an integral scheme and $\mathcal{J}$ be a coherent $\Ss{X}$-submodule of $\mathcal{R}_X$ different from zero. If all the positive powers of $\mathcal{J}$ are complete, then the scheme obtained by blowing-up $\mathcal{J}$ is normal.
\end{lemma}

\begin{lemma}[{\cite[Lemma~5.3]{Li}}]
\label{lema:induccioncomplete}
Let $X$, $\mathcal{R}_X$ and $\mathcal{J}$ be as in Lemma~\ref{lema:blowupcomplete}. Let
\begin{equation*}
 f\colon X \to Y,   
\end{equation*}
be a quasi-compact, quasi-separated birational morphism. If $\mathcal{J}$ is complete, then so is $f_* \mathcal{J}$.
\end{lemma}

\section{Adapted resolutions}
\label{Sec:Adap}
In this section we introduce the notion of adapted resolutions of maximal Cohen-Macaulay modules over any dimension. We prove that if the singularity has dimension two, then the minimal adapted resolution is an adapted resolution (Proposition~\ref{Prop:2daconstruccion}).

\begin{definition}
Let $(X,x)$ be the complex analytic germ of a normal $n$-dimensional singularity. Let $M$ be a maximal Cohen-Macaulay $\Ss{X}$-module. A resolution 
\begin{equation*}
    \pi \colon \Rs \to X,
\end{equation*}
is called an \emph{adapted resolution} associated to $M$ if $\pi^* M / \mathrm{tor}$ is a locally free $\Ss{\Rs}$-module.
\end{definition}

Note that, if $(X,x)$ is a normal surface singularity and $M$ is a reflexive $\Ss{X}$-module, then its minimal adapted resolution is an adapted resolution. The following proposition tell us that adapted resolutions exist in any dimension.

\begin{proposition}
\label{prop:BlowUpM2}
Let $(X,x)$ be the germ of a normal $n$-dimensional singularity. Let $M$ be a maximal Cohen-Macaulay module and 
\begin{equation*}
    f\colon \mathrm{Bl}_M(X) \to X,
\end{equation*}
be the blow-up of $X$ at $M$. Then, any resolution 
\begin{equation*}
    \sigma\colon Y \to \mathrm{Bl}_M(X),
\end{equation*}
is an adapted resolution associated to $M$.
\end{proposition}
\proof
The module $M$ is a maximal Cohen-Macaulay, therefore $M$ is free over the regular part of $X$. Hence, the morphism $f$ is an isomorphism over the regular part of $X$. Thus
\begin{equation*}
    f\circ\sigma\colon Y \to X,
\end{equation*}
is a resolution of $X$. Denote by
\begin{equation*}
 \tilde{M}:= f^*M / \mathrm{tor} \quad \text{ and } \quad  \tilde{\Sf{M}} := \sigma^* \tilde{M},
\end{equation*}
 where tor denotes the torsion part of $f^*M$. 
 
 The sheaf $\tilde{M}$ is locally free and generated by global sections, therefore $\tilde{\Sf{M}}$ is also locally free and generated by global sections. Now, consider the following two exact sequences:
\begin{align}
\label{exact:Prop3.2Ex1}
    &0 \to \mathrm{tor} \to f^* M \to \tilde{M} \to 0,\\
\label{exact:Prop3.2Ex2}
    &0 \to \mathrm{tor} \to \sigma^*f^* M \to \sigma^* f^* M / \mathrm{tor} \to 0.
\end{align}
Applying the functor $\sigma^*$ to the exact sequence~\eqref{exact:Prop3.2Ex1} and using~\eqref{exact:Prop3.2Ex2}, we get the following commutative diagram:
\begin{equation*}
    \begin{tikzcd}
        0 \arrow[r] & \ker \arrow[r] \arrow[dr, phantom, "\circlearrowleft"]& \sigma^* f^* M \arrow[r] \arrow[dr, phantom, "\circlearrowleft"] & \tilde{\Sf{M}} \arrow[r] & 0\\
        0 \arrow[r] & \mathrm{tor} \arrow[u] \arrow[r] & \sigma^* f^* M \arrow[r] \arrow[u,"="] & \sigma^* f^* M / \mathrm{tor} \arrow[r]\arrow[u,"\alpha"] & 0
    \end{tikzcd}
\end{equation*}
Notice that the morphism $\alpha$ is an isomorphism. Indeed, it is clearly a surjection. Now, its kernel is torsion. But $\sigma^* f^* M / \mathrm{tor}$ is torsion-free. Therefore, $\alpha$ is also injective. Since $\alpha$ is an isomorphism, the sheaf $\sigma^* f^* M / \mathrm{tor}$ is locally free. This proves the proposition.
\endproof

As an application of Theorem~\ref{Th:BlowUpM} we can construct the minimal adapted resolution of a reflexive module as follows:

\begin{proposition}
\label{Prop:2daconstruccion}
Let $(X,x)$ be the germ of a normal surface singularity. Let $M$ be a reflexive $\Ss{X}$-module and 
\begin{equation*}
    \pi\colon \Rs \to X,
\end{equation*}
be the associated minimal adapted resolution. Let $f \colon Bl_M(X) \to X$ be the blow-up of $X$ at the module $M$ and 
\begin{equation*}
    \rho\colon \widetilde{\mathrm{Bl}_M(X)}_{\text{min}} \to \mathrm{Bl}_M(X),
\end{equation*}
be the minimal resolution of $\mathrm{Bl}_M(X)$. Then, $\widetilde{\mathrm{Bl}_M(X)}_{\text{min}} \cong \Rs$.
\end{proposition}
\proof
Consider the following commutative diagram:
\begin{equation*}
\begin{tikzpicture}
  \matrix (m)[matrix of math nodes,
    nodes in empty cells,text height=2ex, text depth=0.25ex,
    column sep=3.5em,row sep=3em] {
    \widetilde{\mathrm{Bl}_M(X)}_{\text{min}} & \Rs\\
    \mathrm{Bl}_M(X) & X\\
};
\draw[-stealth] (m-2-1) edge node[below]{$f$} (m-2-2);
\draw[-stealth] (m-1-1) edge node[left]{$\rho$} (m-2-1);
\draw[-stealth] (m-1-2) edge node[right]{$\pi$} (m-2-2);
\draw[-stealth] (m-1-2) edge node[auto]{$\phi$} (m-2-1);
\draw[-stealth] (m-1-2) edge node[above]{$\varphi$} (m-1-1);
\end{tikzpicture}
\end{equation*}
where $\phi$ is given by the universal property of the blow-up of $X$ at $M$ and $\varphi$ comes from the universal property of the minimal resolution of $\mathrm{Bl}_M(X)$. By~\cite[Proposition~5.1]{BoRo}, the resolution $\Rs$ is the minimal resolution of $X$ such that the full sheaf $\Sf{M}:= \left(\pi^* M\right)^{\smvee \smvee}$ is generated by global sections. By Proposition~\ref{prop:BlowUpM2}, the full sheaf $\tilde{\Sf{M}}:=  \left(\rho^* f^* M\right)^{\smvee \smvee}$ is generated by global sections, hence the morphism $\varphi$ is an isomorphism. This proves the proposition.
\endproof

\section{Normality of the blow-up}
\label{sec:Blowup}
In the previous section we used the blow-up of a reflexive module to recover its minimal adapted resolution, so it is natural to study the properties of this blow-up. In this section we prove that the Raynaud-Gruson blow-up of a special module over a Gorenstein, normal, surface singularity is a partial resolution (Theorem~\ref{Th:BlMXNormal}). Furthermore, we prove that this blow-up can be constructed via a resolution (depending on the reflexive module) and the first Chern class of the associated full sheaf (Theorem~\ref{Theo:BlowUpMMinAdapt}). 

First, we observe that the pushforward via a resolution map commutes with tensor product modulo torsion of certain coherent modules.

\begin{lemma}
\label{lemma:DtensorM}
Let $(X,x)$ be the germ of a normal surface singularity. Let $\pi \colon \Rs \to X$ be any resolution with $E$ the exceptional divisor. Let $\Sf{M}$ be a $\Ss{\Rs}$-module generated by global sections and $\Sf{A}$ be an $\Ss{\Rs}$-module of dimension one such that its support intersects $E$ in a finite number of points. Then, the natural morphism
\begin{equation}
    \alpha \colon \pi_* \Sf{M} \otimes \pi_* \Sf{A} \to \pi_* \left( \Sf{M} \otimes \Sf{A}\right),
\end{equation}
is a surjection.
\end{lemma}

\proof
Consider the natural morphism
\begin{equation*} 
\alpha \colon \pi_* \Sf{M} \otimes \pi_* \Sf{A} \to \pi_* \left(\Sf{M}\otimes \Sf{A}\right).
\end{equation*}
Since the support of $\Sf{A}$ intersects the exceptional divisor in a finite number of points, we can identify $\pi_* \left(\Sf{M}\otimes \Sf{A}\right)$ with $\Sf{M}\otimes \Sf{A}$. 

Let  $m \otimes a$ be a section of $\Sf{M}\otimes \Sf{A}$. Since $\Sf{M}$ is generated by global sections there exist global sections $\psi_1, \dots, \psi_n$ of $\Sf{M}$ and sections $f_1, \dots, f_n$ of $\Ss{\Rs}$ defined in an open neighbourhood of the support of $\Sf{A}$ such that $m = \sum_i \psi_i f_i$. Denote by $m' \otimes a' = \sum (\psi_i \otimes f_i\cdot a)$.
By construction we get that $\alpha(m' \otimes a')$ is $m \otimes a$. Therefore, the natural morphism $\alpha$ is a surjection.
\endproof

From now on, let $(X,x)$ be the  germ of a normal Gorenstein surface singularity and $M$  be a special module of rank $r$. Let $\pi \colon \Rs \to X$ be the minimal adapted resolution associated  to $M$ with $E=\cup E_i$ the exceptional divisor with irreducible components $E_i$. Let
\begin{equation*}
 \Sf{M}:= \left(\pi^* M \right)^{\smvee \smvee} \quad \text{and} \quad \Sf{L}:=\mathrm{det}(\Sf{M}).
\end{equation*} By~\cite[Proposition~7.4]{BoRo}, the full sheaf $\Sf{M}$ is an extension of the determinant bundle $\Sf{L}$ by $\Ss{\Rs}^{r-1}$.

Take $r$ generic global sections $\phi_1,...,\phi_r$ of $\Sf{M}$ and consider the following exact sequence given by the
sections:
\begin{equation*}
0 \to \Ss{\Rs}^r \xrightarrow{(\phi_1,...,\phi_r)} \Sf{M} \to \Sf{A}' \to 0.
\end{equation*}
By \cite[Lemma~5.4]{BoRo}, the degeneracy module $\Sf{A'}$ is isomorphic to $\Ss{D}$, where $D\subset\Rs$ is a smooth curve meeting the exceptional divisor transversely at its smooth locus.

The following lemma tell us that the norm of $M$ has a representative given by the global sections of $\Sf{L}:=\mathrm{det}(\Sf{M})$.
\begin{lemma}
\label{lema:NormaML}
Let $(X,x)$ be the germ of a normal Gorenstein surface singularity and $M$  be a special $\Ss{X}$-module. Let $\pi \colon \Rs \to X$ be the minimal adapted resolution associated to $M$. Denote by $\Sf{M}=(\pi^*M)^{\smvee \smvee}$ and by $\Sf{L}:=\mathrm{det}(\Sf{M})$. Then, $\pi_* \Sf{L}$ is a representative of $\|M\|_{\Ss{X}}$.
\end{lemma}
\proof
The proof is the same as \cite[Lemma~4.1]{Gus}. 
\endproof

We want to prove that the blow-up of a special module is a partial resolution. By Lemma~\ref{lema:NormaML}, Remark~\ref{remark:BlowupM}, Lemma~\ref{lema:blowupcomplete} and Lemma~\ref{lema:induccioncomplete}, it is enough to prove that $(\pi_* \Sf{L})^n$  is complete for any positive integer $n$ (where $(-)^n$  denotes the product of fractional ideals). This is going to be the strategy of the proof of Theorem~\ref{Th:BlMXNormal}.

The following remark tell us that the tensor product of locally free sheaves coincides with the product as a fractional ideal sheaves. We need this remark in order to prove that  $(\pi_* \Sf{L})^n$  is complete.
\begin{remark}
\label{remark:Ltensorn}
Let $(X,x)$ be the germ of a normal singularity. Let $\pi \colon \Rs \to X$ be a resolution and let $\Sf{L}$ be a line bundle over $\Rs$. Notice that the natural morphism
\begin{equation}
\label{eq:remark1}
    \Sf{L}^{\otimes n} \to \Sf{L}^{n},
\end{equation}
is an isomorphism. Indeed, the morphism~\eqref{eq:remark1} is always a surjection. Now, the sheaf $\Sf{L}$ is locally free, hence it is flat. Therefore, the natural morphism is injective.
\endproof
\end{remark}

The following theorem tell us that the blow-up of a special module is a partial resolution.

\begin{theorem}
\label{Th:BlMXNormal}
Let $(X,x)$ be the germ of a normal Gorenstein surface singularity and $M$  be a  special module. Let $f \colon \mathrm{Bl}_M(X) \to X$ be the blow-up of $X$ at the module $M$. Then, $\mathrm{Bl}_M(X)$ is normal.
\end{theorem}
\proof
Let
\begin{equation*}
 \pi \colon \Rs \to X,   
\end{equation*}
be the minimal adapted resolution associated to $M$. Denote by $\Sf{M}=(\pi^*M)^{\smvee \smvee}$ and by $\Sf{L}$ its determinant. The resolution is small with respect to the Gorenstein form. Hence, the canonical cycle $Z_K$ is non-negative. Moreover, $D$ does not meet the support of $Z_K$ (see \cite[Proposition~5.14]{BoRo}). Therefore, for any positive integer $n$ we have  
\begin{equation*}
\text{Tor}_1^{\Ss{\Rs}}(\Ss{D}, \mathcal{L}^{\otimes n} \otimes  \Ss{Z_K}) = 0 \quad \text{and} \quad \Ss{D} \otimes \left( \mathcal{L}^{\otimes n} \otimes  \Ss{Z_K}\right) = 0.    
\end{equation*}
By these equalities, applying $- \otimes \left( \Sf{L}^{\otimes n} \otimes \Ss{Z_K}\right)$ to the exact sequence
\begin{equation}
\label{exctseq:detM1}
0 \to \Ss{\Rs} \to \Sf{L} \to \Ss{D} \to 0,
\end{equation}
we get
\begin{equation}
\label{eq:isosDet2}
\Sf{L}^{\otimes n} \otimes \Ss{Z_K} \cong \Sf{L}^{\otimes (n+1)} \otimes \Ss{Z_K}.
\end{equation}
Tensoring \eqref{exctseq:detM1}  with $\Sf{L}^{\otimes n}$ we get
\begin{equation}
\label{exctseq:detM2}
0 \to \Sf{L}^{\otimes n} \to \Sf{L}^{\otimes (n+1)} \to \Ss{D}\otimes \Sf{L}^{\otimes n} \to 0.
\end{equation}
Now, applying the functor $\pi_*$ to the exact sequence \eqref{exctseq:detM2} and using the isomorphism \eqref{eq:isosDet2} we obtain the following commutative diagram:
\begin{equation}
\label{diagram:Th4.31}
\begin{tikzcd}
\dots \arrow[r] & \pi_* \left(\Ss{D}\otimes \Sf{L}^{\otimes n}\right) \arrow[r] & R^1 \pi_* \Sf{L}^{\otimes n} \arrow[r] \arrow[d] \arrow[dr, phantom, "\circlearrowleft"] &R^1 \pi_* \Sf{L}^{\otimes (n+1)} \arrow[r] \arrow[d] & 0 \\
&0 \arrow[r]&  R^1 \pi_* \left(\Ss{Z_K}\otimes \Sf{L}^{\otimes n}\right) \arrow[r] \arrow[d] & R^1 \pi_* \left(\Sf{L}^{\otimes (n+1)} \otimes \Ss{Z_K}\right) \arrow[r] \arrow[d] & 0\\
&  &  0 & 0 & 
\end{tikzcd}
\end{equation}
where the two columns are induced by the exact sequence
\begin{equation}
\label{eq:ExactCanonico}
    0 \to \Cs{\Rs} \to \Ss{\Rs} \to \Ss{Z_K} \to 0.
\end{equation}
Note that the existence of the exact sequence~\eqref{eq:ExactCanonico} follows from the fact that the the resolution is small with respect to the Gorenstein form.

By the diagram~\eqref{diagram:Th4.31} and induction on $n$ (the case $n=1$ follows by~\cite[Lemma~7.5]{BoRo}) we get the following equalities:
\begin{equation}
\label{eq:Th4.3.1}
    \dimc{R^1 \pi_* \left(\Sf{L}^{\otimes n}\right)} = \dimc{R^1 \pi_* \left( \Ss{Z_K}\otimes \Sf{L}^{\otimes n}\right)} =p_g.
\end{equation}
Now, applying the functor $\pi_*$ to the exact sequence~\eqref{exctseq:detM1} and by~\eqref{eq:Th4.3.1}  we get the exact sequence:
\begin{equation}
\label{exctseq:abajodetM1'}
0 \to \Ss{X} \to \pi_* \Sf{L} \to \pi_* \Ss{D} \to 0.
\end{equation}
The exact sequences \eqref{exctseq:detM1} and \eqref{exctseq:abajodetM1'} tell us that $\Sf{L}$ is generated by global sections. Thus, the sheaf $\Sf{L}^{\otimes n}$ is also generated by global sections. Now, we prove by induction that the natural morphism
\begin{equation}
\label{eq:Th4.3.2}
    \left(\pi_* \Sf{L} \right)^{\otimes n} \to \pi_* \left( \Sf{L}^{\otimes n} \right),
\end{equation}
is a surjection for any positive integer $n$. The case $n=1$ is clearly true. Assume that the assertion is true for some $n=k$. We then prove the assertion for $n=k+1$. Consider the following commutative diagram obtained by tensoring \eqref{exctseq:detM1} and \eqref{exctseq:abajodetM1'} with $\Sf{L}^{\otimes k}$ and $\left(\pi_* \Sf{L}\right)^{\otimes k}$, respectively:
\begin{equation*}
\begin{tikzcd}
0 \arrow[r] & \pi_* \left( \Sf{L}^{\otimes k}\right) \arrow[r] & \pi_* \left( \Sf{L}^{\otimes (k+1)}\right) \arrow[r,"\sigma"] & \pi_* \left( \Ss{D} \otimes \Sf{L}^{\otimes k} \right)\arrow[r] & 0\\
& \left(\pi_*  \Sf{L}\right)^{\otimes k} \arrow[r] \arrow[u,"\alpha_1"] \arrow[ur, phantom, "\circlearrowleft"]  & \left(\pi_*  \Sf{L}\right)^{\otimes (k+1)} \arrow[r] \arrow[u,"\alpha_2"] \arrow[ur, phantom, "\circlearrowleft"] & \pi_* \Ss{D} \otimes \left(\pi_*  \Sf{L}\right)^{\otimes k} \arrow[r] \arrow[u,"\alpha_3"]& 0 
\end{tikzcd}
\end{equation*}
where the morphisms in the columns are the natural ones. The morphism $\sigma$ is a surjection by~\eqref{eq:Th4.3.1}. The morphism $\alpha_1$ is a surjection by the induction hypothesis. The morphism $\alpha_3$ is also a surjection by Lemma~\ref{lemma:DtensorM}. This implies that $\alpha_2$ is a surjection. This proves the induction step. This proves the surjectivty of~\eqref{eq:Th4.3.2}. 

Now, we prove that $(\pi_* \Sf{L})^n$  is complete for any positive integer $n$. By Remark~\ref{remark:Ltensorn}, we have 
\begin{equation}
\label{eq:Ln}
    \pi_* \left( \Sf{L}^{\otimes n}\right) = \pi_* \left(\Sf{L}^n \right).
\end{equation}
Consider the composition
\begin{equation}
\label{eq:composition}
\begin{tikzpicture}
  \matrix (m)[matrix of math nodes,
    nodes in empty cells,text height=1.5ex, text depth=0.25ex,
    column sep=2.5em,row sep=2em] {
  \left(\pi_* \Sf{L} \right)^{\otimes n} & \pi_* \left( \Sf{L}^{\otimes n} \right) & \pi_* \left( \Sf{L}^{n} \right). \\
};
\draw[-stealth] (m-1-1) -- (m-1-2);
\draw[-stealth] (m-1-2) edge node[auto]{$\sigma$} (m-1-3);
\end{tikzpicture}
\end{equation}
Since the natural map given in~\eqref{eq:Th4.3.2} is a surjection and by the equality~\eqref{eq:Ln}, we get that the composition~\eqref{eq:composition} is also a surjection. Therefore, $(\pi_* \Sf{L})^n = \pi_* \left(\Sf{L}^n \right)$. By \eqref{eq:Ln}, this implies $\pi_* \left( \Sf{L}^{\otimes n} \right) = \left(\pi_* \Sf{L} \right)^n$. Finally by Lemma~\ref{lema:induccioncomplete} the ideal $\left(\pi_* \Sf{L} \right)^n$ is complete. Now, the theorem follows by Remark~\ref{remark:BlowupM}, Lemma~\ref{lema:blowupcomplete} and Lemma~\ref{lema:induccioncomplete}. This proves the theorem.
\endproof

The following theorem tell us how to recover the normalization of the blow-up at a reflexive module using the minimal adapted resolution.
 
\begin{theorem}
\label{Theo:BlowUpMMinAdapt}
 Let $(X,x)$ be the germ of a normal surface singularity. Let $M$ be a reflexive  $\Ss{X}$-module of rank $r$. Let $\pi\colon \Rs \to X$ be the minimal adapted resolution associated to $M$ with exceptional divisor $E$. Let $E_1,\dots,E_n$ be the irreducible components of $E$ and $\Sf{M}:=(\pi^*M)^{\smvee \smvee}$ be the full sheaf associated to $M$. Then, the normalization of $\mathrm{Bl}_M(X)$ is obtained by contracting the irreducible components $E_j$ such that $c_1\left(\Sf{M} \right) \cdot E_j =0$.
\end{theorem}
\proof
Let  $E_1,\dots,E_m$ be the irreducible components of $E$ such that  $c_1(\Sf{M})\cdot E_j \neq 0$. Let 
\begin{equation*}
    h\colon \Rs \to Y,
\end{equation*}
be the contraction of all the irreducible components of $E$ different from $E_1,\dots,E_m$. Denote by $E':=\bigcup_{j=1}^m E_j$ and $S:=h(E\setminus E')$. Notice that $S$ is a finite set of cardinality equal to the number of connected components of $E\setminus E'$. By Grauert's contraction theorem~\cite{Gr} the variety $Y$ is normal. Let 
\begin{equation*}
    v \colon Y \to X,
\end{equation*} be the natural morphism such that $\pi = v \circ{} h$. 

We prove that $h_* \Sf{M}$ is a locally free $\Ss{Y}$-module generated by its global sections. Let $\phi_1, \dots, \phi_r$ be generic global sections of $\Sf{M}$. By~\cite[Lemma~5.4]{BoRo}, the exact sequence given by the sections is 
\begin{equation}
\label{eq:Th4.8.1}
    0 \to \Ss{\Rs}^r \xrightarrow{(\phi_1,...,\phi_r)} \Sf{M} \to \Ss{D} \to 0.
\end{equation}
Applying the functor $h_*$ to the exact sequence~\eqref{eq:Th4.8.1} we get
\begin{equation} \label{eq:hM}
    0 \to \Ss{Y}^r \to h_*\Sf{M} \stackrel{\psi}{\longrightarrow} h_*\Ss{D}.
\end{equation}
Notice that the morphism $h$ is an isomorphism in the complement of $E \setminus E'$ and the support of $D$ only intersects $E'$. Hence, the support of $h_*\Ss{D}$ and $S$ are disjoint sets. Thus, the morphism $\psi$ is a surjection and $h_* \Sf{M}$ is locally free sheaf generated by its global sections. Now, denote by  $\tilde{\Sf{M}}:= \left(v^* M \right)^{\smvee \smvee}$ and consider the natural morphism
\begin{equation*}
    \theta \colon h_* \Sf{M} \to \tilde{\Sf{M}}.
\end{equation*}
The kernel and cokernel of $\theta$ are torsion, $Y$ is normal and $h_* \Sf{M}$ is locally free. Thus, the morphism $\theta$ is injective. We now prove that $\theta$ is also a surjection. Consider the following exact sequence:
\begin{equation*}
    0 \to h_* \Sf{M} \stackrel{\theta}{\longrightarrow} \tilde{\Sf{M}} \to K \to 0.
\end{equation*}
Applying the functor $\Homs_{\Ss{Y}}\left(-, \Cs{Y} \right)$ to the exact sequence we get
\begin{equation*}
    0\to \Homs_{\Ss{Y}}\left(\tilde{\Sf{M}}, \Cs{Y} \right)   \to \Homs_{\Ss{Y}}\left(h_* \Sf{M}, \Cs{Y} \right) \to \Exts_{\Ss{Y}}^1 \left(K, \Cs{Y} \right) \to 0.
\end{equation*}
The sheaf $K$ is supported in a finite set. Therefore,  by~\cite[Theorem~3.3.10]{BrHe} the sheaf $\Exts_{\Ss{Y}}^1 \left(K, \Cs{Y}\right)$ must vanish. Thus, we get
\begin{equation}
\label{eq:OmegaisoY}
    \Homs_{\Ss{Y}}\left(\tilde{\Sf{M}}, \Cs{Y} \right) \cong  \Homs_{\Ss{Y}}\left(h_* \Sf{M}, \Cs{Y} \right).
\end{equation}
Since $h_* \Sf{M}$ is locally free and $\tilde{\Sf{M}}$ is reflexive. We conclude that both sheaves are $\Cs{Y}$-reflexive modules. Recall that $\Cs{Y}$-reflexivity is equivalent to reflexivity. Therefore,
\begin{equation*}
    h_* \Sf{M} \cong \tilde{\Sf{M}}.
\end{equation*}
Consequently, the sheaf $\tilde{\Sf{M}}$ is locally free and it is generated by its global sections. Since $\tilde{\Sf{M}}$ is generated by global sections, we get
\begin{equation*}
    \tilde{\Sf{M}} \cong v^* M / \mathrm{tor}.
\end{equation*}
In particular, the sheaf $v^* M / \mathrm{tor}$ is locally free. Let
\begin{equation*}
 n \colon \mathrm{NBl}_M(X) \to \mathrm{Bl}_M(X),
\end{equation*}
be the normalization of $\mathrm{Bl}_M(X)$. Denote by $ \widetilde{\mathrm{NBl}_M(X)}_{\text{min}}$ the minimal resolution of  $\mathrm{NBl}_M(X)$. By the universal property of the blow-up of $X$ at $M$ and the universal property of the normalization, there exist morphisms $g\colon Y \to \mathrm{Bl}_M(X)$ and $\gamma \colon Y \to \mathrm{NBl}_M(X)$ such that the diagram commutes
\begin{equation*}
\begin{tikzpicture}
  \matrix (m)[matrix of math nodes,
    nodes in empty cells,text height=1ex, text depth=0.35ex,
    column sep=3.5em,row sep=3em] {
    \widetilde{\mathrm{NBl}_M(X)}_{\text{min}}& \mathrm{NBl}_M(X) & \mathrm{Bl}_M(X)\\
    \Rs & Y &\\
};
\draw[-stealth] (m-1-1) edge node[auto]{} (m-1-2);
\draw[-stealth] (m-2-1) edge node[auto]{$h$} (m-2-2);
\draw[-stealth] (m-2-2) edge node[auto]{$g$} (m-1-3);
\draw[-stealth] (m-2-2) edge node[auto]{$\gamma$} (m-1-2);
\draw[-stealth] (m-1-2) edge node[auto]{$n$} (m-1-3);
\end{tikzpicture}
\end{equation*}
By Proposition~\ref{Prop:2daconstruccion} we know that $\Rs \cong \widetilde{\mathrm{NBl}_M(X)}_{\text{min}}$. Hence, the morphism $\gamma$ is an isomorphism. This proves the theorem.
\endproof

\begin{corollary}
\label{cor:BlowUpMMinAdapt}
Let $(X,x)$ be the germ of a normal Gorenstein surface singularity. Let $M$ be a special $\Ss{X}$-module. Then, the blow-up of $X$ at $M$ has at worst normal Gorenstein singularities.
\end{corollary}
\proof
Let
\begin{equation*}
 \pi\colon \Rs \to X,   
\end{equation*}
be the minimal adapted resolution associated to $M$ with exceptional divisor $E$ and $\Sf{M}:=(\pi^*M)^{\smvee \smvee}$ be the full sheaf associated to $M$. Let  $E_1,\dots,E_m$ be the irreducible components of $E$ such that $c_1(\Sf{M})\cdot E_j \neq 0$. If the singularity is rational, then the corollary follows by~\cite[Theorem~2.2]{Curto}. Now assume that the singularity is not rational. Let 
\begin{equation*}
    h\colon \Rs \to Y,
\end{equation*}
be the contraction of all the irreducible components of $E$ different from $E_1,\dots,E_m$. Denote by $E':=\bigcup_{j=1}^m E_j$ and $S:=h(E\setminus E')$. Recall that $S$ is a finite set.

By Theorem~\ref{Th:BlMXNormal} and Theorem~\ref{Theo:BlowUpMMinAdapt} we get $Y \cong \mathrm{Bl}_M(X)$. By Grauert's conraction theorem~\cite{Gr} the restriction
\begin{equation*}
h|_{\Rs \setminus (E\setminus E')}\colon \Rs \setminus (E\setminus E') \to Y\setminus S,    
\end{equation*}
is a isomorphism. By~\cite[Remark~5.3]{BoRo}, the Gorenstein $2$-form $\Omega_{\Rs}$ does not have any zero or pole along $\Rs \setminus (E\setminus E')$. Therefore, there exists a two form that does no vanishes over $Y\setminus S$, i.e., $Y$ has only Gorenstein singularities. This proves the corollary.
\endproof

\begin{remark}
\label{Remark:Unico}
Let $(X,x)$ be the germ of a normal Gorenstein surface singularity. Let $M$ be a special $\Ss{X}$-module. By Corollary~\ref{cor:BlowUpMMinAdapt}, the blow-up of $X$ at $M$ has at worst normal Gorenstein singularities. The only case when the blow-up of $X$ at $M$ is smooth is in the case of the $A_1$-singularity. In all other cases the blow-up of $X$ at $M$ always has normal Gorenstein singularities. 
\end{remark}

\begin{corollary}
\label{cor:BlowUpMMinAdapt2}
Let $(X,x)$ be the germ of a normal Gorenstein surface singularity. Then, two special modules give isomorphic partial resolutions if and only if they are isomorphic up to free summands.
\end{corollary}
\proof
The proof follows by Theorem~\ref{Theo:BlowUpMMinAdapt}, Proposition~ \ref{Prop:2daconstruccion} and \cite[Corollary~7.12]{BoRo}.
\endproof

Now, by Theorem~\ref{Th:BlMXNormal} and Theorem~\ref{Theo:BlowUpMMinAdapt} we generalize the McKay correspondence as follows: 

\begin{corollary}
\label{Cor:PrincipalNuevo}
Let $(X,x)$ be a normal Gorenstein surface singularity. Then, there exists a bijection between the following sets: 
\begin{enumerate}
\item The set of special, indecomposable $\Ss{X}$-modules up to isomorphism.
\item The set of irreducible divisors $E$ over $x$, such at any resolution of $X$ where $E$ appears, the Gorenstein form has neither zeros nor poles along $E$.
\item The set of partial resolutions $\psi\colon Y \to X$ with irreducible exceptional divisor $E$ such that:
    \begin{enumerate}
        \item the partial resolution is dominated by a resolution which is small with respect to the Gorenstein form.
        \item the Gorenstein form does not have any zeros or poles along $E \setminus \text{Sing $Y$}$.
    \end{enumerate}
\end{enumerate}
\end{corollary}
\proof
Clearly, if $M_1$ and $M_2$ are two special modules such that $M_1 \cong M_2$, then the blow-up of $X$ at $M_1$ and at $M_2$ are isomorphic. By Corollary~\ref{cor:BlowUpMMinAdapt}, the blow-up of $X$ at $M_1$ is a Gorenstein partial resolution with irreducible exceptional divisor. Hence (1) implies (3).

Let $\psi\colon Y \to X$ be a partial resolution satisfying the properties of (a) and (b). Therefore, there exists a resolution small with respect to the Gorenstein form such that the Gorenstein form does not have zeros or poles along $E'$, where $E'$ is the strict transform of $E$. Hence (3) implies by (2). The bijection between (1) and (2) follows from \cite[Corollary~7.12]{BoRo}. This proves the corollary.
\endproof

\section{Applications via matrix factorizations}
\label{Sec:Matrix}
In this section we use the theory of matrix factorizations  in order to compute examples of blow-ups at reflexive modules. We recall some preliminaries on matrix factorizations.

\subsection{Matrix factorizations}
From now on, let $(S,\mathfrak{m})$ be a regular local ring and suppose that $R=S/I$ is Henselian, where $I$ is a principal ideal of $S$ generated by $f$.

Recall, the notion of matrix factorization given by Eisenbud~\cite{Ei} gives an equivalence of categories between maximal Cohen-Macaulay $R$-modules and pairs of matrices with entries in $S$ satisfying some conditions. We review this construction. See~\cite{Ei} or~\cite[Chapter~7]{Yos} for more details.

\begin{definition}
Let $R$ be the coordinate ring of the a hypersurface defined by $f$ in $(S,\mathfrak{m})$. A \emph{matrix factorization} of $f$ is an ordered pair of $n \times n$-matrices $(\Phi, \Psi)$ with entries in $S$ such that 
\begin{equation*}
    \Phi \cdot \Psi = f \cdot \mathrm{Id}_{S^n}, \quad \Psi \cdot \Phi = f \cdot \mathrm{Id}_{S^n},
\end{equation*}
where $\mathrm{Id}_{S^n}$ is the $n\times n$ identity matrix.

A \emph{morphism} between matrix factorizations $(\Phi_1, \Psi_1)$ and $(\Phi_2, \Psi_2)$ is a pair of $n \times n$ matrices $(\alpha, \beta)$ with entries in $S$ such that
\begin{equation*}
    \alpha \cdot \Phi_1 = \Phi_2 \cdot \beta, \quad \beta \cdot \Psi_1 = \Psi_2 \cdot \alpha.
\end{equation*}
A matrix factorization is \emph{reduced} if and only if
\begin{equation*}
    \im \Phi \subset \mathfrak{m} S^n \quad \text{and} \quad \im \Psi \subset \mathfrak{m} S^{n}.
\end{equation*}
\end{definition}

Using matrix factorizations Eisenbud~\cite{Ei} proved the following:
\begin{theorem}[{\cite{Ei}}]
\label{Th:EisTh}
There is a one-to-one correspondence between:
\begin{enumerate}
    \item equivalence classes of reduced matrix factorizations of $f$.
    \item isomorphism classes of non-trivial periodic minimal free resolutions of $R$-modules of periodicity two.
    \item maximal Cohen-Macaulay $R$-modules without free summands.
    \end{enumerate}
\end{theorem}
We sketch the idea of a part of the proof of Theorem~\ref{Th:EisTh} which will be used later in this section. Let $M$ be a maximal Cohen-Macaulay $R$-module without free summands. By the Auslander-Buchsbaum-Serre theorem we have
\begin{equation*}
    \text{proj dim}_S M = \dim S - \text{depth } M = 1. 
\end{equation*}
Thus, there is a free resolution of $M$ as $S$-module of length $1$:
\begin{equation*}
    0 \to S^n \stackrel{\Phi}{\longrightarrow} S^n \to M \to 0.
\end{equation*}
Since $f$ annihilates the module $M$, we get $f \cdot S^n \subset \im \Phi$. Hence, there exists a matrix $\Psi$ with entries in $S$ such that 
\begin{equation*}
    \Phi \cdot \Psi = f \cdot \mathrm{Id}_{S^n},
\end{equation*}
where the pair $(\Phi, \Psi)$ is a $n\times n$-matrix factorization of $f$ with $\coker{\Phi} =M$.

Conversely, let $(\Phi, \Psi)$ be a matrix factorization of $f$. Denoting $\bar{\Phi}$ and $\bar{\Psi}$ the matrices $\Phi, \Psi$ modulo $(f)$ respectively, we have the following exact sequence of $R$-modules:
\begin{equation*}
    \dots \to R^n \stackrel{\bar{\Psi}}{\longrightarrow} R^n \stackrel{\bar{\Phi}}{\longrightarrow} R^n \stackrel{\bar{\Psi}}{\longrightarrow}  \dots
\end{equation*}
If $f$ is not a zero-divisor in $S$, then the complex
\begin{equation*}
    \dots \to R^n \stackrel{\bar{\Psi}}{\longrightarrow} R^n \stackrel{\bar{\Phi}}{\longrightarrow} R^n \to   \coker{\bar{\Phi}} \to 0,
\end{equation*}
is exact. Hence it is a periodic free resolution of $\coker{\bar{\Phi}}$ with periodicity two. Furthermore, the module $\coker{\bar{\Phi}}$ is a maximal Cohen-Macaulay $R$-module.

\subsection{The blow-up at a matrix factorization.}
In this subsection we use the blow-up given by Villamayor~\cite{Villa} and the matrix factorizations in the case of $S=\mathbb{C}\{x,y,z\}$ and $R$ the coordinate ring of a normal hypersurface, i.e., $R=S/(f)$ with $f \in S$ and $R$ is a normal ring. Recall that any hypersurface singularity is a Gorenstein singularity. 

Let $M$ be a reflexive $R$-module of rank $r$. By~\cite[Proposition~2.5]{Villa} we can use the matrix factorization associated to $M$ to obtain a representative of $\|M\|_R$. Some computations can be done by hand, but in some cases we use the software {\sc Singular}~{4-1-2}~\cite{SING} and the libraries resolve.lib~\cite{RES} and sing.lib~\cite{SingLib}.

\begin{example}
\label{ex:Ex1}
Let $f=xy+z^{n+1}$, i.e., $R$ is the $A_n$-singularity. The matrix factorization of rational double points are well known, see for example \cite{Kaj}. In the case of the $A_n$-singularity the matrix factorizations are
\begin{equation}
\Phi_k
=
\begin{bmatrix}
    y & -z^{n+1-k} \\
    z^k & x
\end{bmatrix}
, \quad
\Psi_k
=
\begin{bmatrix}
    x & z^{n+1-k}\\
    -z^k & y \\
\end{bmatrix},
\end{equation}
where $k$ is a integer such that $0\leq k \leq n$. Using the matrix $\Phi_k$ we get the following morphism:
\begin{equation*}
    R^2 \stackrel{\Phi_k}{\longrightarrow} R^2 \to  M(\Phi_k) \to 0,
\end{equation*}
where $M(\Phi_k):=\coker{\Phi_k}$. Let $K$ be the kernel of the morphism $\Phi_k$ and denote by
\begin{equation*}
    K_1 := \set{\begin{bmatrix}
    y\\
    z^k
\end{bmatrix}
\cdot g \in R^2}{g \in R}.
\end{equation*}
By~\cite[Proposition~2.5]{Villa}, the ideal 
\begin{equation*}
    I_k=(y,z^k),
\end{equation*}
is a representative of $\|M(\Phi_k)\|_R$. Therefore, the blow-up at $M(\Phi_k)$ of $R$ is the blow-up at the ideal $I_k$. In this case $\mathrm{Bl}_{M(\Phi_k)}(R)$ has at most two singular points:
\begin{enumerate}
    \item If $\mathrm{Bl}_{M(\Phi_k)}(R)$ has one singular point, then the singularity is $xy+z^n$.
    \item If $\mathrm{Bl}_{M(\Phi_k)}(R)$ has two singular points, then one singularity is $xy+z^{n-l}$ and the other singularity is $xy+z^{l+1}$ with $1 \leq l \leq n-2$.
\end{enumerate}
Both cases were exactly as predicted by Corollary~\ref{Cor:PrincipalNuevo} or \cite{Gus}.
\end{example}
We now consider a different singularity. From now on, let $f=x^3 + y^3 + z^3$ and $R=\mathbb{C}\{x,y,z\}/(f)$. Hence $R$ is a normal Gorenstein surface singularity. In this case all the reflexive modules were classified by Kahn~\cite{Ka} and all the special modules were classified in \cite{BoRo}. Several people have studied this singularity and the category of reflexive modules, for example \cite{Ka,Laza1,Laza2}. First, we study the blow-up at the fundamental module of $R$.

\begin{definition}[{\cite[Definition~11.5]{Yos}}]
The \emph{fundamental exact sequence} of $R$ is the following exact sequence (unique, up to non-canonical isomorphism):
\begin{align*}
    0 \to R \to A \to R \to \mathbb{C} \to 0,
\end{align*}
corresponding to a non-zero element of $\Ext_R^2\left(R/\mathfrak{m},R\right) \cong \mathbb{C}$. The module $A$ is called the \emph{fundamental module} of $R$.
\end{definition}

The fundamental module of $R$ is an indecomposable reflexive module  of rank $2$ (see \cite[Chapter~11]{Yos} for more properties about the fundamental module). A natural question is the following: Is the fundamental module, a special module? We show:

\begin{proposition}
Let $f=x^3 + y^3 + z^3$ and $R=\mathbb{C}\{x,y,z\}/(f)$. Then, the fundamental module of $R$ is not special.
\end{proposition}
\proof
The idea is the same as in Example~\ref{ex:Ex1}: we use the matrix factorization of $A$ to compute the ideal that we need to blow-up. Then, we use Corollary~\ref{cor:BlowUpMMinAdapt} in order to check if the module is special. The matrix factorization associated to $A$ was computed by Yoshino and Kawamoto in~\cite{Yos1} and Laza, Pfister and Popescu in~\cite{Laza1}. The periodic free resolution of $A$ is the following:
\begin{equation}
    \dots \to R^4 \stackrel{\Psi_A}{\longrightarrow} R^4 \stackrel{\Phi_A}{\longrightarrow} R^4 \to A \to 0,
\end{equation}
where
\begin{equation}
\Phi_A
=
\begin{bmatrix}
    x^2 & -y  & -z   & 0 \\
    y^2 & x   & 0    & -z \\
    z^2 & 0   & x    & y \\
    0   & z^2 & -y^2 & x^2 
\end{bmatrix}
, \quad
\Psi_A
=
\begin{bmatrix}
    x & y & z & 0 \\
    -y^2 & x^2 & 0 & z \\
    -z^2 & 0 & x^2 & -y \\
    0 & -z^2 & y^2 & x 
\end{bmatrix}.
\end{equation}
By~\cite[Proposition~2.5]{Villa}, we can choose as a representative of $\|A\|_R$ the ideal generated by all the $2\times2$-minors of the matrix given by
\begin{equation}
D
=
\begin{bmatrix}
    -z   & 0 \\
    0    & -z \\
    x    & y \\
    -y^2 & x^2 
\end{bmatrix}.
\end{equation}
By Theorem~\ref{Th:BlowUpM}, the blow-up of $R$ at $A$ is the blow-up at the ideal
\begin{equation}
\label{eq:IdealAuslander}
    I= (z^2,yz,xz,x^3+y^3).
\end{equation}
Using {\sc Singular}~{4-1-2}~\cite{SING} and the libraries resolve.lib~\cite{RES} and sing.lib~\cite{SingLib} one can check that the blow-up at the ideal $I$ has $4$ smooth charts, therefore the blow-up of $R$ at $A$ is a resolution. Then, by Corollary~\ref{cor:BlowUpMMinAdapt} and Remark~\ref{Remark:Unico} the module $A$ is not special.
\endproof

We now give one example where the blow-up is not normal.
\begin{example}
\label{ex:Ex2}
Consider the following matrix factorization of $f$:
\begin{equation}
\Phi
=
\begin{bmatrix}
    x+ y & -z^2 \\
    z & x^2-xy+y^2 
\end{bmatrix}
, \quad
\Psi
=
\begin{bmatrix}
     x^2-xy+y^2 & z^2\\
     -z& x+y  \\
\end{bmatrix}
\end{equation}
Let $M=\coker{\Phi}$. Then, using {\sc Singular}~{4-1-2}~\cite{SING} and the libraries resolve.lib~\cite{RES} and sing.lib~\cite{SingLib} one can check that the blow-up of $R$ at $M$ does not have an isolated singularity, hence it is not normal.
\end{example}

\subsection*{Acknowledgments}
We would like to thank Javier Fern\'andez de Bobadilla and Ananyo Dan for helpful discussions during the course of this work.

\end{document}